\documentclass[11pt]{article}
\usepackage{mathrsfs}

\newcounter{example}

\newtheorem{theorem}{Theorem}[section]

\usepackage{latexsym}
\usepackage{amsmath}
\usepackage{amssymb}
\usepackage{graphicx}
\usepackage{textcomp}
\usepackage{multirow}
\usepackage{enumerate}
\usepackage{cancel}
\usepackage{hyperref}
\usepackage[usenames,dvipsnames]{color}
\usepackage[left, pagewise]{lineno}
\usepackage{epstopdf}
\usepackage{algpseudocode}
\usepackage{algorithm}

\title{Controlling an Alien Predator Population by Regional Controls}
\author{Sebastian Ani\c{t}a\thanks{Corresponding author. Faculty of Mathematics, 
``Alexandru Ioan Cuza'' University of Ia\c{s}i, and 
``Octav Mayer'' Institute of Mathematics of the Romanian Academy, Ia\c{s}i 700506, Romania. Email: {sanita@uaic.ro}.}
\and Vincenzo Capasso\thanks{ADAMSS (Centre for Advanced Applied Mathematical and Statistical Sciences), 
Universit\'a degli Studi di Milano, 20133 Milano, Italy.  Email: {vincenzo.capasso@unimi.it}.}
\and Gabriel Dimitriu\thanks{Department of Medical Informatics and Biostatistics,
University of Medicine and Pharmacy ``Grigore T. Popa'',
Ia\c{s}i 700115, Romania. Email: {gabriel.dimitriu@umfiasi.ro}. }
}
\date{}
\begin{document}
\maketitle

\begin{abstract}
We investigate the problem of minimizing the total cost of the
damages produced by an alien predator population and of the
regional control paid to reduce  this population. The dynamics of
the predators is described by a prey-predator system with either
local or nonlocal reaction terms. A sufficient condition for the
zero-stabilizability (eradicability) of predators is given in
terms of the sign of the principal  eigenvalue of an appropriate
operator that is not self-adjoint, and a stabilizing feedback
control with a very simple structure is indicated. The
minimization related to such a feedback control is treated for a
closely related minimization problem viewed as a regional control
problem. The level set method is a  key  ingredient. An iterative
algorithm to decrease the total cost is obtained and numerical
results show the effectiveness of the theoretical results.

A spatially structured SIR problem may be described by the same
system; in this case  the above mentioned minimization problem is
related to the problem of eradication of  an epidemic by regional
controls.
\end{abstract}

{\bf Keywords:} Zero-stabilization; regional control; prey-predator system; SIR system.

\section{Setting of the problem}\label{AD_s1}
Consider the following reaction--diffusion system which describes the dynamics of two interacting populations: prey and predator that are
free to move in the habitat $\Omega\subset\mathbb{R}^2$ and are subject to a control acting in a subset $\omega\subset\Omega$ on the predators.
\begin{equation}\label{AD_e1}
\left\{
\begin {array}{ll}
\partial _t h(x,t)-d_1\Delta h(x,t)=r(x)h(x,t)-\rho(x)h(x,t)^2 - h(x,t)(Bp(\cdot,t))(x), &  (x,t)\in Q, \\
\partial _t p(x,t)-d_2\Delta p(x,t)=-a(x)p(x,t)  + c_0h(x,t)(Bp(\cdot,t))(x) + \chi_{\omega }(x)u(x,t), & (x,t)\in Q, \\
\partial _{\nu}h(x,t)=\partial _{\nu}p(x,t)=0, & (x,t)\in \Sigma, \\
h(x,0)=h_0(x),\quad p(x,0)=p_0(x), & x\in \Omega\,.
\end{array}
\right.
\end{equation}
Here $\Omega$ is a bounded domain (open and connected) with a sufficiently smooth boundary $\partial\Omega$, $\omega$ is an open subset,
$Q=\Omega \times (0,+\infty)$, $\Sigma=\partial\Omega \times (0, + \infty)$.
$h(x,t)$ and $p(x,t)$ are the spatial densities at time $t$ of the prey, respectively predator populations. The diffusion coefficients $d_1$ and $d_2$
are positive constants, $r(x)$ is the growth rate and $\rho(x)h(x,t)^2$ is a local logistic term of preys. $a(x)$ is the decreasing rate of the
predator population. The quantity $h(x,t)(Bp(\cdot,t))(x)$ gives the density of captured prey population at position $x$, which is transformed
 into biomass via a conversion rate $c_0\in(0,+\infty)$. Possible choices of the operator $B\in L(L^2(\Omega))$  will  be  discussed later.

The homogeneous Neumann boundary conditions describe the  no flux
of populations across  the boundary of the habitat.  $h_0(x)$ and
$p_0(x)$ denote the initial prey and predator population densities
at position $x,$  respectively. The control $u$ acts on the
predators only, in the  subregion $\omega$; $\chi_\omega$ is the
characteristic function of $\omega$.

We assume that
\begin{itemize}
\item[{\bf (A1)}] $r, \rho, a, h_0, p_0\in L^{\infty }(\Omega)$, $r(x)\geq r_0,\ \rho(x)\geq\rho_0$ a.e. $x\in\Omega$
($r_0$, $\rho_0$ are positive constants);
$$h_0(x)\geq 0, \quad p_0(x)\geq 0 \qquad \mbox{\rm a.e. }x\in \Omega ,$$
and $h_0$ and $p_0$ are not identically zero.
\item[{\bf (A2)}] $B\in L(L^2(\Omega ))\cap L(L^{\infty }(\Omega ))$, $(By)(x) \geq 0$ a.e. $x\in\Omega $, for any $y\in L^2(\Omega)$ such that $y(x)\geq 0$ a.e. $x\in\Omega $.
\end{itemize}
\vspace{2mm}

Two  cases are of particular interest to us:
\vspace{3mm}

CASE $1$. If $(By)(x)=c(x)y(x)$ for $y\in L^2(\Omega)$, where
$c\in L^\infty(\Omega), \ c(x)\geq 0$ a.e. $x\in\Omega$, then the
functional  response to predation is of the  usual Lotka-Volterra
type. \vspace{3mm}

CASE $2$.   If $(By)(x)=\int_{\Omega }\kappa(x,x')y(x')\,dx'$ for
$y\in L^2(\Omega )$, where $\kappa \in L^{\infty }(\Omega \times
\Omega)$, $\kappa (x,x')\geq 0$ a.e. $(x,x')\in\Omega \times
\Omega$, then the functional  response to predation is such that
predators, however coming from any position $x',$ upon  predation
at position $x$ will stay and produce offsprings at this new
position (the predators follow the prey). For other prey-predator systems with nonlocal terms see \cite{gva}.\vspace{3mm}

We  may notice that system (\ref{AD_e1}) may also model  a
spatially structured SIR epidemic system, in which case  $h(x,t)$
and $p(x,t)$ represent the spatial density of the susceptible and
infective population, respectively. With  $c_0=1,$ CASE $1$
presented above  corresponds to a local infection rate, while CASE
$2$   corresponds to a  nonlocal infection rate as proposed by
D.G. Kendall \cite{kendall} (see also \cite{aronson}, and
\cite{capasso}).

For the SIR system $r(x)=b(x)-\mu(x)$, where $b(x)$ is the birth
rate and $\mu(x)$ is the natural death rate at position $x$;
$a(x)=\mu(x)+\tilde{\mu}(x)$, where $\tilde{\mu}(x)$ is the
additional removal rate due  to the extra death rate caused by the
disease and the possible  natural recovery  rate. The
infective population do not have offsprings and the recovered
individuals acquire  immunity. The control $u$ may describe the
additional removal of infectives (because of either  recovery  by
treatment or isolation)  due to a planned regional intervention by
the relevant public health system.

If we view $p$ as an alien pest population density or as an
infective population density, then it is of great interest to know
if there exists a control $u$ such that, for   the solution
$(h^u,p^u)$ to (\ref{AD_e1}),   $\lim_{t\to\infty}p^u(\cdot,t)=0$
in an appropriate functional space. \vspace{2mm}

\bf Definition. \rm The predator population is {\it zero-stabilizable} (eradicable),
if for any $h_0$, $p_0$ satisfying (A1), there exists $u\in L^\infty_{loc}(\overline{\omega}\times [0,+\infty))$ such that
\begin{equation}\label{AD_e2}
h^u(x,t),\ p^u(x,t)\ge 0\quad \mbox{\rm a.e. } (x,t)\in Q\,,
\end{equation}
and
\begin{equation}\label{AD_e3}
\lim_{t\to +\infty}p^u(\cdot,t)=0\quad \mbox{\rm in } L^\infty(\Omega).
\end{equation}
We are dealing with zero-stabilizability (eradicability) with state constraints.
\vspace{3mm}

We will see that the problem of eradicability is deeply related to the sign of the principal eigenvalue $\lambda^\omega_{1\gamma}$ for
\begin{equation}\label{AD_e4}
\left\{
\begin {array}{ll}
-d_2\Delta\Psi(x)+a(x)\Psi(x)+\gamma\chi_\omega(x)\Psi(x)-c_0K(x)(B\Psi)(x)=\lambda\Psi(x), &  x\in\Omega, \\
\partial _{\nu}\Psi(x)=0, & x\in \partial\Omega,
\end{array}
\right.
\end{equation}
where $\gamma\in [0, +\infty)$, and $K$ is the unique maximal nonnegative solution to
\begin{equation}\label{AD_e5}
\left\{
\begin {array}{ll}
-d_1\Delta K(x) =r(x) K(x) - \rho(x) K(x)^2,  &  x\in\Omega\,, \\
\partial _{\nu} K(x)=0, & x\in \partial\Omega\,.
\end{array}
\right.
\end{equation}
(actually (\ref{AD_e5}) has two nonnegative solutions, the trivial one and $K$; see \cite{AFL}).

Notice that the elliptic operator in (\ref{AD_e4}) is not
self-adjoint, so we cannot use the same arguments as in
\cite{AK17} to derive the basic properties of
$\lambda^\omega_{1\gamma}$ and of the corresponding eigenspace. We
may  however apply  the Krein-Rutman Theorem  (for details see the
Appendix).

Our strategy will be to diminish the pest (resp. infective)
population using a bilinear control with a very simple structure
$u:=-\gamma p$. Here $\gamma\in [0, +\infty)$ represents a
constant affordable predator  elimination (resp. treatment  or
isolation) rate. If we consider the feedback control $u:=-\gamma
p$\ ($\gamma\in [0,+\infty$)), then Banach's fixed point theorem
implies that (\ref{AD_e1}) has a unique solution $(h,p)$ which has
nonnegative components.

The following zero-stabilizability (eradicability) shall be proved in the next section and extends a result in \cite{AFL} .
\begin{theorem}\label{AD_t1}
If $\lambda^\omega_{1\gamma}>0$, where $\gamma\in [0,+\infty)$, then the feedback control $u:=-\gamma p$ realizes {\rm(\ref{AD_e2})} and {\rm(\ref{AD_e3})}, for any $h_0$, $p_0$ satisfying {\rm(A1)}.
\end{theorem}

We shall see that  moreover,
$$\lim_{t\to +\infty} p(\cdot,t)=0\quad \mbox{\rm  in } L^\infty(\Omega)\,,$$
at the rate of $e^{-\lambda^\omega_{1\gamma}t}\,.$

From the  point of view of Geometric Measure Theory, in   $2D,$
the geometry of $\omega$ can be described  by its three Minkowski
functionals. Actually in this work we shall consider only the area
and the perimeter of $\omega$ (denoted  by
$length(\partial(\omega))$). We assume that the cost to be paid in
order to act in $\omega$ is
$$\alpha\cdot area(\omega) + \beta\cdot length(\partial\omega)\,,$$
where $\alpha, \beta > 0$, and is paid once  for all  for
installing the harvesting devices (resp. for treatment or
isolation units) (see also \cite{AK17}).

After all our goal is to find a  subregion $\omega$ which
minimizes the total cost of the damages produced by the pest
population  (resp. the cost of the treatment for the infective
population),  and of the costs associated with  the intervention
in $\omega$, namely
$$ Minimize \ \left \{ \theta \int_0^{\infty} \int_{\Omega} h(x,t) (Bp(\cdot,t))(x)\,dx\,dt+\alpha \cdot area (\omega)+\beta \cdot length (\partial \omega) \right \}, \leqno{\bf (\tilde{P})}$$
subject to $\omega$, where $\theta $ is a positive constant.

We shall see that usually $h(x,t)\leq K(x)$ a.e. in $Q$, and $p(x,t)\le y(x,t)$ a.e. in $Q$, where $y$ is the solution to
\begin{equation}\label{AD_e6}
\left\{
\begin {array}{ll}
\partial _t y(x,t)-d_2\Delta y(x,t)= -a(x)y(x,t) + c_0 K(x)(By(\cdot,t))(x) \\
\qquad\qquad\qquad\qquad\qquad\ \  - \gamma \chi_\omega(x)y(x,t)\,, & (x,t)\in Q\,, \\
\partial _{\nu} y(x,t)=0, & (x,t)\in \Sigma\,,\\
y(x,0)=p_0(x), & x\in \Omega\,.
\end{array}
\right.
\end{equation}
This implies that
$$\theta \int_0^{\infty} \int_{\Omega} h(x,t) (Bp(\cdot,t))(x)\,dx\,dt
\leq \theta \int_0^{\infty} \int_{\Omega} K(x)(By(\cdot,t))(x)\,dx\,dt\,,$$
and consequently, if for a certain $\omega$ we get a small value of the cost functional in
$$ Minimize \ \left \{ \theta \int_0^{\infty} \int_{\Omega} K(x) (By(\cdot,t))(x)\,dx\,dt+\alpha \cdot area (\omega)+\beta \cdot length (\partial \omega) \right \}, \leqno{\bf (P)}$$
(subject to $\omega$), then for the same $\omega$ we get even a smaller value for the cost functional in ${\bf(\tilde{P})}$.

We shall treat problem ({\bf P}) as a shape optimization problem,
using the level set method (see \cite{osher_fedkiw} and references
therein). We shall use the implicit interface according to which
$\partial\omega$ is the zero-isocontour of a certain function
$\varphi:\overline{\omega}\to\mathbb{R}$ and $\omega=\{ x\in
\Omega; \ \varphi (x)>0\}$, while $\partial \omega =\{ x\in
\overline{\Omega}; \varphi (x)=0\}$. If $\varphi $ is an implicit
function of $\omega$, then
$$ area(\omega) = \int_\Omega H(\varphi(x))\,dx\, \quad \mbox{\rm and } \quad length (\partial \omega)=\int_{\Omega}\delta (\varphi (x))|\nabla \varphi (x)|\,dx.$$
Here $H$ is the Heaviside function and $\delta$ is its  Dirac
Delta generalized derivative. Hence, we may rewrite problem ({\bf
P}) as
$$ Minimize \left \{ \  \theta \int_0^{\infty}\int_{\Omega}K(x)(By(\cdot ,t))(x)\,dx\,dt +
\alpha\int_{\Omega }H(\varphi (x))\,dx +\beta \int_{\Omega}\delta
(\varphi (x))|\nabla \varphi (x)|\,dx \right \}, \leqno{\bf (P)}$$
subject to $\varphi $, where $y$ is the solution to
\begin{equation*}
\left\{
\begin {array}{ll}
\partial _t y(x,t)-d_2\Delta y(x,t)= -a(x)y(x,t) + c_0 K(x)(By(\cdot,t))(x)- \gamma H(\varphi(x))y(x,t)\,, & (x,t)\in Q\,, \\
\partial _{\nu} y(x,t)=0, & (x,t)\in \Sigma\,,\\
y(x,0)=p_0(x), & x\in \Omega\,.
\end{array}
\right.
\end{equation*}

For stabilization problems related to reaction-diffusion systems in biology we refer to \cite{anita_VK_1}-\cite{anita_VK_8}. Some optimal control problems in mathematical
biology have been treated in \cite{anita}, \cite{anita2}, \cite{AD_2008}, \cite{ADS_2014}, \cite{bv}, \cite{bressan_coclite}, \cite{coclite_measure}, \cite{coclite_garavello_spinolo}, 
\cite{fl}, \cite{GDS_2010}, \cite{hy}, \cite{lenhart_toronto}-\cite{llw}, \cite{SD_2012}.
We have also to mention some recent results concerning the regional control in population dynamics \cite{AK15}.
For basic notions and methods in shape optimization theory see \cite{bucur}, \cite{chanv}, \cite{delfour_zolesio}, \cite{getreuer}, \cite{henrot}, \cite{osher_fedkiw}, \cite{sethian}, \cite{sokolowski}.

\vspace{2mm}

The present paper contains  a continuation of  the investigations
started in \cite{AK17}, in  the more complicated  situation  of a
non-selfadjoint operator $B$. We shall use   here  most of the
notations adopted  in \cite{AK17}.

Here is the plan of the paper. Section 2 concerns the proof of
Theorem~\ref{AD_t1}. Section 3 is devoted to an ``approximation''
of problem ({\bf P}). An iterative algorithm to decrease the total
cost by changing $\omega$ is derived in Section 4. Some numerical
results for a nonlocal interaction are given. Final remarks are
presented in the next section. Basic properties of
$\lambda^\omega_{1\gamma}$ and of the corresponding eigenfunctions
are proved in the Appendix.

\section{Proof of Theorem~\ref{AD_t1}}
Let $\gamma\in [0,+\infty)$ such that $\lambda^\omega_{1\gamma}>0$.
For any sufficiently small $\varepsilon>0$ we have that $\lambda^\omega_{1\gamma}(\varepsilon)>0$, where $\lambda^\omega_{1\gamma}(\varepsilon)$ is the principal eigenfunction to (\ref{AD_e7})
(we have used that $\lim_{\varepsilon\rightarrow 0} \lambda^\omega_{1\gamma}(\varepsilon) = \lambda^\omega_{1\gamma}$; see the Appendix).
\begin{equation}\label{AD_e7}
\left\{
\begin {array}{ll}
-d_2\Delta\Psi(x)+a(x)\Psi(x)+\gamma\chi_\omega(x)\Psi(x)-c_0(K(x)+\varepsilon)(B\Psi )(x)=\lambda\Psi(x), &  x\in\Omega, \\
\partial _{\nu}\Psi(x)=0, & x\in \partial\Omega\,.
\end{array}
\right.
\end{equation}
Consider $\Psi_1$ an eigenfunction to (\ref{AD_e7}) corresponding to the eigenvalue $\lambda^\omega_{1\gamma}(\varepsilon)$ and satisfying
$$\Psi_1(x)>0\,,\quad\forall x\in\overline{\Omega}$$
(notice that $\Psi_1\in C(\overline{\Omega})$). It follows that there exists $\zeta\in(0,+\infty)$ such that
$$\Psi_1(x) \geq \zeta,\quad\forall x\in\overline{\Omega}\,.$$
Since the solution $(h,p)$ to (\ref{AD_e1}), corresponding to $u:=-\gamma p$, satisfies
$$0\leq h(x,t)\leq \tilde{h}(x,t)\quad \mbox{\rm a.e.\ in}\ Q\,,$$
where $\tilde{h}$ is the unique solution to
\begin{equation*}
\left\{
\begin {array}{ll}
\partial _t h(x,t)-d_1\Delta h(x,t)=r(x)h(x,t)-\rho(x)h(x,t)^2\,, &  (x,t)\in Q, \\
\partial _{\nu}h(x,t) = 0\,, & (x,t)\in \Sigma, \\
h(x,0)=h_0(x)\,, & x\in \Omega\,.
\end{array}
\right.
\end{equation*}
(this follows in a standard manner using the fact that $0\le h(x,t)(Bp(\cdot,t))(x)$ a.e. in $Q$; for comparison results for parabolic equations see \cite{barbu1}, \cite{friedman}, \cite{protter}), and using that

$$\lim_{t\to\infty}\tilde{h}(\cdot,t)=K \quad\mbox{in}\ L^\infty(\Omega)\,,$$
we may conclude that for any sufficiently small $\varepsilon>0$, there exists $T(\varepsilon)\in (0,+\infty)$ such that
$$h(x,t)\leq K(x)+\varepsilon\quad\mbox{\rm a.e. in }\Omega \times (T(\varepsilon ),+\infty )\, ,$$
and consequently
$$0\leq c_0 h(x,t)(B\Psi)(x)\le c_0(K(x)+\varepsilon)(B\Psi)(x)\quad \mbox{\rm a.e. in }\Omega \times (T(\varepsilon ),+\infty ) \,,$$
for any $\Psi\in L^\infty(\Omega)$, $\Psi (x)\ge 0$ a.e. $x\in\Omega $.
Since $p(\cdot,T(\varepsilon))\in L^\infty(\Omega)$, it follows that there exists $\tau\in (0,+\infty)$ such that
$$0\leq p(x,T(\varepsilon))\leq \tau\Psi_1(x)\quad \mbox{\rm a.e. } x\in\Omega\,.$$

On the other hand, the function $z(x,t)=\tau\Psi_1(x)e^{-\lambda_{1\gamma}^{\omega }(\varepsilon)t}$ is the solution to
\begin{equation*}
\left\{
\begin {array}{ll}
\partial _t z(x,t)-d_2\Delta z(x,t)= - a(x) z(x,t) + c_0 (K(x) + \varepsilon) (Bz(\cdot,t))(x) \\
\qquad\qquad\qquad\qquad\qquad\ \ - \gamma \chi_\omega(x) z(x,t)\,, &  (x,t)\in \Omega\times (T(\varepsilon), +\infty)\,, \\
\partial _{\nu}z(x,t) = 0\,, & (x,t)\in \partial\Omega\times (T(\varepsilon), + \infty)\,,\\
z(x,T(\varepsilon))=\tau\Psi_1(x)\,, & x\in\Omega\,,
\end{array}
\right.
\end{equation*}
and in a standard manner we may conclude the following comparison result:
$$0\leq p(x,t)\leq \tau \Psi_1(x) e^{-\lambda^\omega_{1\gamma}(\varepsilon)t}\quad\mbox{a.e.\ in}\ \Omega\times (T(\varepsilon ), + \infty)\,.$$
This implies that $\lim_{t\to +\infty} p(\cdot, t)=0$ in $L^\infty(\Omega)$ at the rate of $e^{-\lambda^\omega_{1\gamma}(\varepsilon)t}\,.$
\vspace{2mm}

\bf Remark. \rm
Usually, even without the presence of predators the initial density $h_0(x)$ is less or equal than $K(x)$ for $x\in\Omega$. In this case we get that
$h(x,t)\leq K(x)$ a.e. in $Q$, and consequently we may infer (this follows in a standard way; see \cite{friedman}, \cite{protter})
$$0\le p(x,t) \le \tilde{p}(x,t)\quad\mbox{a.e.\ in}\ Q\,,$$
where $\tilde{p}$ is the solution to
\begin{equation*}
\left\{
\begin {array}{ll}
\partial _t \tilde{p}(x,t)-d_2\Delta \tilde{p}(x,t)= - a(x) \tilde{p}(x,t) + c_0 K(x) (B\tilde{p}(\cdot,t))(x)- \gamma \chi_\omega(x) \tilde{p}(x,t)\,, &  (x,t)\in Q\,, \\
\partial _{\nu}\tilde{p}(x,t) = 0\,, & (x,t)\in \Sigma, \\
\tilde{p}(x,0)=\tau_0 \tilde{\Psi}_1(x)\,, & x\in \Omega\,,
\end{array}
\right.
\end{equation*}
$\tau_0\in (0, +\infty)$ and $\tilde{\Psi}_1$ is an eigenfunction to (\ref{AD_e4}) corresponding to the eigenvalue $\lambda^\omega_{1\gamma}$, satisfying
$$\tilde{\Psi}_1(x)\ge\zeta_0\,,\quad\forall x\in\overline{\Omega}\,.$$
Here $\zeta_0\in (0, +\infty)$\ (actually, $\tilde{\Psi}_1\in C(\overline{\Omega})$), and $p_0(x)\leq \tau _0\zeta _0$ a.e. $x\in \Omega $.
We may infer that
$$0\le p(x,t)\le \tau_0\tilde{\Psi}_1(x) e^{-\lambda^\omega_{1\gamma}t}\quad\mbox{a.e.\ in}\ Q$$
(because $\tilde{p}(x,t)=\tau_0\tilde{\Psi}_1(x) e^{-\lambda^\omega_{1\gamma}t}\,,\ (x,t)\in\overline{Q}$) and that
$\lim_{t\to + \infty} p(\cdot,t)=0$ in $L^\infty(\Omega)$ at the rate of $e^{-\lambda^\omega_{1\gamma}t}$\,.

\section{Regional control of the pest population}\label{AD_s3}

In the sequel we shall denote by  $T>0$  a large number,
$Q_T=\Omega\times (0,T)$, $\Sigma_T=\partial\Omega\times (0,T)$,
$d=d_2$, $y_0=p_0$.
$H_{\sigma}(s)=\frac{1}{2}\left(1+\frac{2}{\pi}\mbox{arctan}\frac{s}{\sigma}\right)$
is a mollified version of $H(s)$ and its derivative
$\delta_\sigma(s)=\frac{\sigma}{\pi(\sigma^2+s^2)}$ is a mollified
version of $\delta(s)$.

 Without loss of generality, we
may assume that $\theta=1$. Assume that $h_0(x)\le K(x)\
\mbox{a.e.}\ x\in\Omega$. Problem (${\bf \tilde{P}}$) may
be ``approximated'' by the following regional control problem:
$$\mbox{Minimize}\ J(\varphi)\,,\leqno{\bf (RC)}$$
where $\varphi:\overline{\Omega}\to\mathbb{R}$ is a smooth function,
\begin{eqnarray*}
J(\varphi) =J_{damage}(\varphi) + \alpha J_{area}(\varphi) + \beta J_{perimeter}(\varphi)\,.
\end{eqnarray*}
Here
$$J_{damage}(\varphi) =  \int_0^T\int_{\Omega}K(x)(By^\varphi(\cdot,t))(x)\,dx\,dt  $$
is the cost of the damages produced by the predators,
$$J_{area}(\varphi )=\int_{\Omega}H_{\sigma }(\varphi (x))\,dx$$
is an approximation of the area of $\omega =\{ x\in \Omega ; \ \varphi(x)>0\}$, and
$$J_{perimeter}(\varphi )=\int_{\Omega}\delta_\sigma(\varphi (x))|\nabla \varphi (x)|\,dx, $$
is an approximation of the perimeter of $\omega $ (length of $\partial \omega $); $y^{\varphi}$ is the solution of
\begin{equation}\label{AD_e8}
\left\{\begin {array}{ll}
\partial _ty(x,t)-d\Delta y(x,t)= -a(x)y(x,t) + c_0 K(x) (By(\cdot,t))(x) - \gamma H_\sigma(\varphi(x))y(x,t)\,, &  (x,t)\in Q_T\,, \\
\partial _{\nu} y(x,t) = 0\,, & (x,t)\in \Sigma_T, \\
y(x,0)=y_0(x)\,, & x\in \Omega\,.
\end{array}
\right.
\end{equation}

\begin{theorem}\label{AD_t2}
For any smooth functions $\varphi, \psi: \overline{\Omega}\to\mathbb{R}$ we have that
\begin{equation}\label{AD_e9}
\begin{array}{ll}
dJ(\varphi )(\psi) & = \displaystyle \int_{\Omega} \delta_{\sigma}(\varphi(x))\psi (x)\left[\gamma \int_0^Tr^{\varphi }(x,t)y^\varphi(x,t)\,dt+\alpha\right. \\
& ~ \\
&\quad -\left.\displaystyle \beta \
\textnormal{div}\left(\frac{\nabla \varphi(x)}{|\nabla
\varphi(x)|}\right)\right]\,dx  + \displaystyle \beta
\int_{\partial \Omega }{{\delta _{\sigma }(\varphi (x))}\over
{|\nabla \varphi (x)|}}\partial _{\nu }\varphi (x)\psi
(x)\,d\ell ,
       \end{array}
\end{equation}
where $r^{\varphi} $ is the solution to
\begin{equation}\label{AD_e10}
\left\{ \begin{array}{ll}
\partial_t r(x,t)+ d\Delta r(x,t)= a(x) r(x,t) - c_0(B^*(K(\cdot)r(\cdot,t)))(x) \\
 \qquad\qquad\qquad\qquad +\gamma H_{\sigma }(\varphi(x))r(x,t)+(B^*K)(x)\,, & (x,t)\in Q_T,\\
\partial _{\nu } r(x,t)=0,  & (x,t)\in\Sigma_T,\\
r(x,T)=0,  & x\in\Omega\,.
           \end{array}
    \right.
\end{equation}
Here $B^*$ is the adjoint of $B\in L(L^2(\Omega))$.
\end{theorem}

{\it Sketch of the proof}. As in the proof of Lemma 3 in
\cite{AK15}, it is possible to prove that, for any smooth
functions $\varphi, \psi :\overline{\Omega }\rightarrow {\bf R},$
we have
$$ \lim_{s\to\ 0}{1\over {s }}[y^{\varphi+s \psi}-y^{\varphi}] = z \quad \mbox{\rm in } C([0,T];L^{\infty }(\Omega )),$$
where $z$ is the solution to the problem
\begin{equation}\label{AD_e11}
\left\{ \begin{array}{ll}
\partial_t z(x,t)- d\Delta z(x,t)= -a(x) z(x,t) + c_0 K(x)(Bz(\cdot ,t))(x) \\
 \qquad\qquad\qquad -\gamma H_{\sigma }(\varphi(x))z(x,t)-\gamma \delta_{\sigma }(\varphi(x)) y^\varphi(x,t)\psi(x)\,, & (x,t)\in Q_T,\\
\partial _{\nu } z(x,t)=0,  & (x,t)\in\Sigma_T,\\
z(x,T)=0,  & x\in\Omega\,.
           \end{array}
    \right.
\end{equation}

If $\varphi , \psi:\overline{\Omega}\to\mathbb{R}$  are arbitrary and smooth functions, then
$$\lim_{s \rightarrow 0}{1\over {s }}[J(\varphi +s\psi )-J(\varphi )]=\int_0^T\int_{\Omega }K(x)(Bz(\cdot,t))(x)\,dx\,dt+  \alpha \int_{\Omega}\delta _{\sigma }(\varphi (x))\psi (x)\,dx$$
$$\begin{aligned}
& + \beta \int_{\Omega} \delta'_{\sigma}(\varphi(x))\psi
(x)|\nabla\varphi(x)|\,dx+\beta \int_{\Omega }\delta _{\sigma
}(\varphi (x)){{\nabla \varphi (x)\cdot \nabla \psi (x)}\over
{|\nabla \varphi (x)|}}\,dx\,.
\end{aligned}$$
After some calculations, we get as in \cite{AK15} that
\begin{equation}\label{AD_e12}
\begin{aligned}
dJ(\varphi )(\psi )=& \int_0^T\int_{\Omega}K(x)(Bz^\varphi(\cdot,t))(x)\,dx\,dt + \alpha \int_{\Omega} \delta_{\sigma}(\varphi(x))\psi(x)\,dx  \\
& - \beta\int_\Omega \delta_\sigma(\varphi(x))\,
\text{div}\left({{\nabla \varphi (x)}\over {|\nabla \varphi
(x)|}}\right)\psi(x)\,dx +\beta \int_{\partial \Omega }{{\delta
_{\sigma }(\varphi (x))}\over {|\nabla \varphi (x)|}}\partial
_{\nu }\varphi (x)\psi (x)\,d\ell \,.
\end{aligned}
\end{equation}
 If we multiply the first equation in (\ref{AD_e10}) by $z$ and integrate over $Q_T$, we obtain after an easy calculation, and using (\ref{AD_e10}) and (\ref{AD_e11}), that
\begin{equation}\label{AD_e13}
\int^T_0\int_\Omega
K(x)(Bz(\cdot,t))(x)\,dx\,dt=\gamma\int^T_0\int_\Omega
\delta_\sigma (\varphi(x))\psi(x) r^\varphi(x,t)
y^\varphi(x,t)\,dx\,dt\,.
\end{equation}
By (\ref{AD_e12}) and (\ref{AD_e13}) we get the conclusion.

\vspace{2mm}

\bf Remark. \rm By (\ref{AD_e9}) we obtain the gradient descent with respect to $\varphi $ (see \cite{getreuer})

\begin{equation}\label{AD_e14}
\left\{
\begin{array}{ll}
      \partial_{s } \varphi(x,s) =  \delta_{\sigma}(\varphi (x,s))\left[ -\gamma \int_0^T r^{\varphi }(x,t)y^\varphi(x,t)\,dt \right. \\
\left. \qquad\qquad\quad -\alpha +\beta \ \text{\rm div}
      \left( \frac{\nabla \varphi(x,s)}{|\nabla \varphi(x,s )|}\right)\right], & x \in \Omega,\ s>0\,, \\
    \displaystyle \frac{\delta_{\sigma}(\varphi(x,s ))}{|\nabla \varphi(x,s )|}\partial_{\nu} \varphi (x,s )=0,  & x\in \partial\Omega,\ s>0\,.
\end{array}
\right.
\end{equation}
\vspace{3mm}

\bf Remark. \rm If $B$ is the operator in CASE 1, then $r^\varphi$ is the solution to
$$\left\{
\begin{array}{ll}
\partial_t r(x,t)+ d\Delta r(x,t)= a(x)r(x,t) - c_0 c(x) K(x) r(x,t) + c(x) K(x) + \gamma H_\sigma (\varphi (x)) r(x,t)\,, & (x,t)\in Q_T,\\
\partial_\nu r(x,t) = 0\,, & (x,t)\in \Sigma_T\,,\\
r(x,T)=0,  & x\in\Omega\,.
\end{array}
\right. $$
\vspace{3mm}

\bf Remark. \rm If $B$ is the operator in CASE 2, then $r^\varphi$ is the solution to
$$\left\{
\begin{array}{ll}
\partial_t r(x,t)+ d\Delta r(x,t)= a(x)r(x,t) - \displaystyle c_0 \int_\Omega \kappa(x',x)K(x')r(x',t)\,dx' \\
\qquad\qquad\qquad + \displaystyle \int_\Omega \kappa(x',x)K(x')\,dx' + \gamma H_{\sigma}(\varphi (x))r(x,t)\,, & (x,t)\in Q_T\,,\\
\partial_\nu r(x,t) = 0\,, & (x,t)\in \Sigma_T\,,\\
r(x,T)=0,  & x\in\Omega\,.
\end{array}
\right. $$

\section{Computational issues}\label{AD_s4}
In this section we approach numerically the optimal control problem for the
model (\ref{AD_e8}) on a two-dimensional domain defining an isolated habitat $\Omega $.\
The systems resulted by discretization were solved iteratively using Matlab's built in function {\sc gmres}, which implements the generalized minimal residual method. 
For our numerical simulations, we have used {\sc gmres} algorithm without ``restarts'' of the iterative method, and found satisfactory to use no preconditioners, and a tolerance for the relative error of $10^{-3}$. 
{\sc gmres} algorithm was also applied to an optimal control problem for a two-prey and one-predator model with diffusion in \cite{ADS_2014}.

The set of the model parameters was asserted to have the
following values:

-- parameters defining the discretization process for space variable
$x=(x_1, x_2)$ in the domain $\Omega=[0,1]\times [0,1]$ and time:
space steps $\Delta x_1=\Delta x_2=2.78$e-$2$, time step $\Delta t=2.78$e-$2$ (36 discretization points on both space and time axes);

-- final time $T$, and maximum number of iterations, $maxiter$: $T=1$, $maxiter = 50$;

-- diffusion parameter: $d=1.$e-$2$;

-- parameters defining birth and mortality rates: $a\equiv 1$, $c_0=1$, $K\equiv 1$, $\gamma=1$;

-- prescribed convergence parameters for $J$ and shape function $\varphi$, given by $\varepsilon_1$ and
$\varepsilon_2$ respectively, and $\sigma$, a parameter for mollified version of
the Heaviside function: $\varepsilon_1=1.$e-$4$,
$\varepsilon_2=1.$e-$5$, and $\sigma=1.$e-$2$;

-- parameters representing the weights in the cost functional $J$: $\theta=1$, and $\alpha$ and $\beta$ vary successively in the set
${\cal{\cal W}}=\{0$,  $0.00001$, $0.0001$, $0.001$, $0.01$, $0.1$, $50$, $75$, $100 \}$.

Theorem \ref{AD_t2} allows us to construct an iterative procedure to update the shape function $\varphi $ defining the subregion where the control acts.
The following three stop criteria have been used in the Algorithm detailed below:
$$iter > maxiter\,,\qquad  \left|J^{(iter+1)} - J^{(iter)}\right| < \varepsilon_1\,,\qquad\mbox{and}\qquad  \|\varphi^{(iter+1)}-\varphi^{(iter)}\|_{L^2(\Omega)}/area(\Omega ) < \varepsilon_2\,.$$

{\small{
\begin{algorithm}[!tbhp]
  \caption{: Iterative scheme to update the shape function $\varphi $ (the subregion where the control acts).}
\label{aprec}
  \begin{algorithmic}[1]
\tt

\State Set $iter := 0$; Choose the positive constants $T$, $d$,
$c_0$, $\gamma$, $maxiter$, $\sigma$ (parameter for mollified
Heaviside function), $\varepsilon_1$ and $\varepsilon_2$.

\State Define the operator $B$, as well as the functions $a$ and $K$.

\State Choose a large value for $J^{(0)}$ and a small constant for $s_0>0$ (artificial time).

\State Initialize shape function: $\varphi^{(0)}:=\varphi^{(0)}(x,0)$\,.

\vspace{1.5mm}

\State Compute  $y^{(iter+1)}$ the solution of \eqref{AD_e8} corresponding to $\varphi^{(iter)}:=\varphi^{(iter)}(x,0)$\,.

\State Compute $\displaystyle J^{(iter+1)}: =\int_0^T\int_{\Omega }K(x)(By^{(iter+1)})(x,t)\,dx\,dt$

\hspace{2.6cm} + $\displaystyle \alpha \int_{\Omega} H_{\sigma
}(\varphi ^{(iter)}(x,0))\,dx$ + $\displaystyle \beta
\int_{\Omega} \delta_{\sigma }(\varphi ^{(iter)}(x,0)) |\nabla
\varphi ^{(iter)}(x,0)|\,dx$.

\State \textbf{IF} $\left|J^{(iter+1)} - J^{(iter)}\right| < \varepsilon_1$
\textbf{THEN} \textbf{STOP}

\hspace{1cm}\textbf{ELSE GO TO} \textbf{8:}

\State Compute $r^{(iter+1)}$ the solution of problem \eqref{AD_e10} corresponding to $\varphi^{(iter)}(\cdot, 0)$ and $y^{(iter+1)}$.

\State Compute $\varphi^{(iter+1)}$ using \eqref{AD_e14} and the initial condition

\hspace{1cm}$\varphi ^{(iter+1)}(x,0):=\varphi ^{(iter)}(x,s_0)$ using a semi-implicit timestep scheme.

\State \textbf{IF} $\|\varphi^{(iter+1)}-\varphi^{(iter)}\|_{L^2(\Omega)}/area(\Omega )<\varepsilon_2$ \textbf{THEN STOP}

\hspace{1cm}\textbf{ELSE} \textbf{IF} $iter > maxiter$ \textbf{THEN STOP}

\hspace{2cm} $iter := iter + 1$

\hspace{2.2cm}\textbf{GO TO} \textbf{5:}

 \end{algorithmic}
\end{algorithm}
}}

\noindent The tolerances $\varepsilon_1 > 0$ in Step 7 and $\varepsilon_2 > 0$ in Step 10 are prescribed convergence parameters. For details about the gradient methods, see \cite{an}.

\bigskip

In what follows, we present results of several numerical simulations corresponding to CASE 2 (Section 1, p. 2), when
$(By)(x)=\int_{\Omega }\kappa(x,x')y(x')\,dx'$ for $y\in L^2(\Omega )$, where $\kappa \in L^{\infty }(\Omega \times \Omega)$, $\kappa (x,x')\geq 0$ a.e. $(x,x')\in\Omega \times \Omega$.
In this case, the numerical response to predation shows that the predators from position $x'$ that captured preys at position $x$ will stay and produce offsprings at this new position (the predators follow the prey).

\bigskip

\textbf{Experiment $1$}.
We have used the following initial levels of the state $y$ and function $\varphi$ for
our numerical simulations:
$y_0(x_1,x_2)=1$, and
$$\varphi_0(x_1,x_2) = \mbox{exp}(-3(x_1-0.5)^2 - 3(x_2-0.5)^2) + \sin(3\pi x_1)\sin(5\pi x_2) - 0.75\,.$$
The function $\kappa(x,x')$ is defined by $\kappa(x,x') \equiv ~|\kappa_1(x_1,x_2) \kappa_2(x'_1,x'_2)|$, where
$$\kappa_1(x_1,x_2) = x_1^2\sin(\pi x_1) + x_2^2\sin(\pi x_2)\,,\qquad\mbox{and}\qquad
\kappa_2(x'_1,x'_2) = 100(x^{\prime 2}_1 \cos(\pi x'_1) +x^{\prime 2}_2\cos(\pi x'_2))\,.$$

\begin{figure}[!tbhp]
\begin{center}
\includegraphics[scale=0.68]{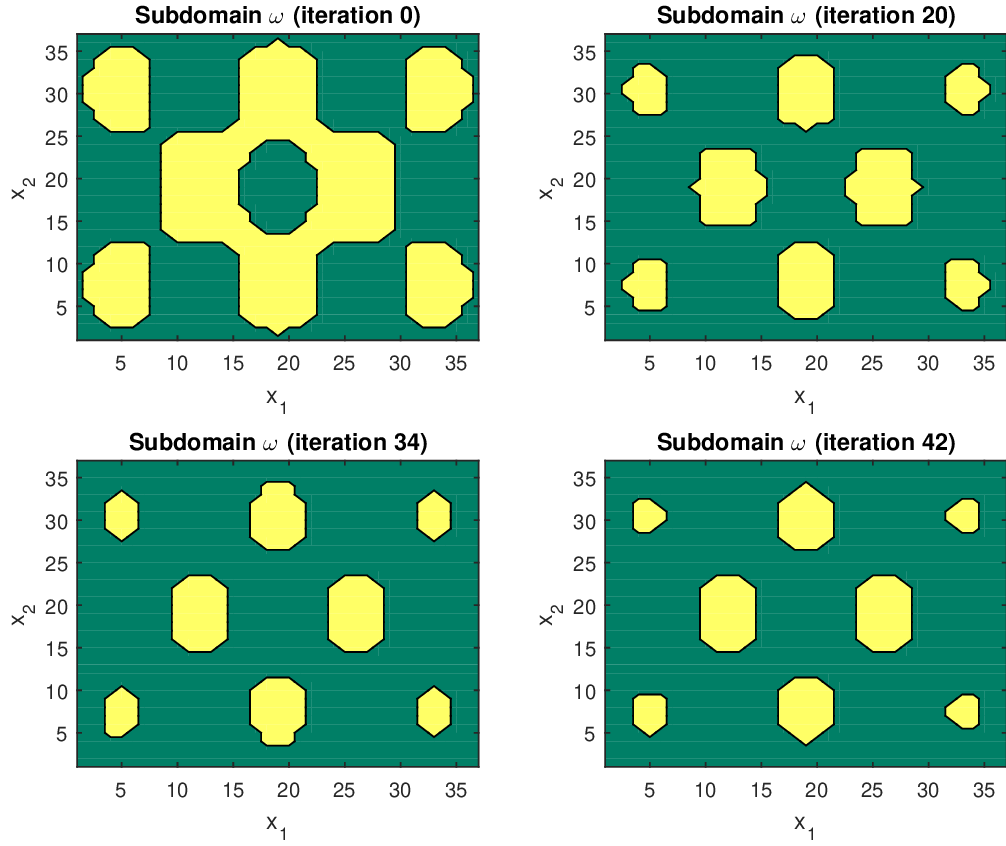}
\caption{Variation of the subdomain  $\omega$ (marked with light color) for $\alpha=100$ and $\beta=0.1$ at iterations $0, 20, 34$ and $42$  (Experiment 1).}
\label{Fig_1}

\bigskip

\bigskip

\includegraphics[scale=0.68]{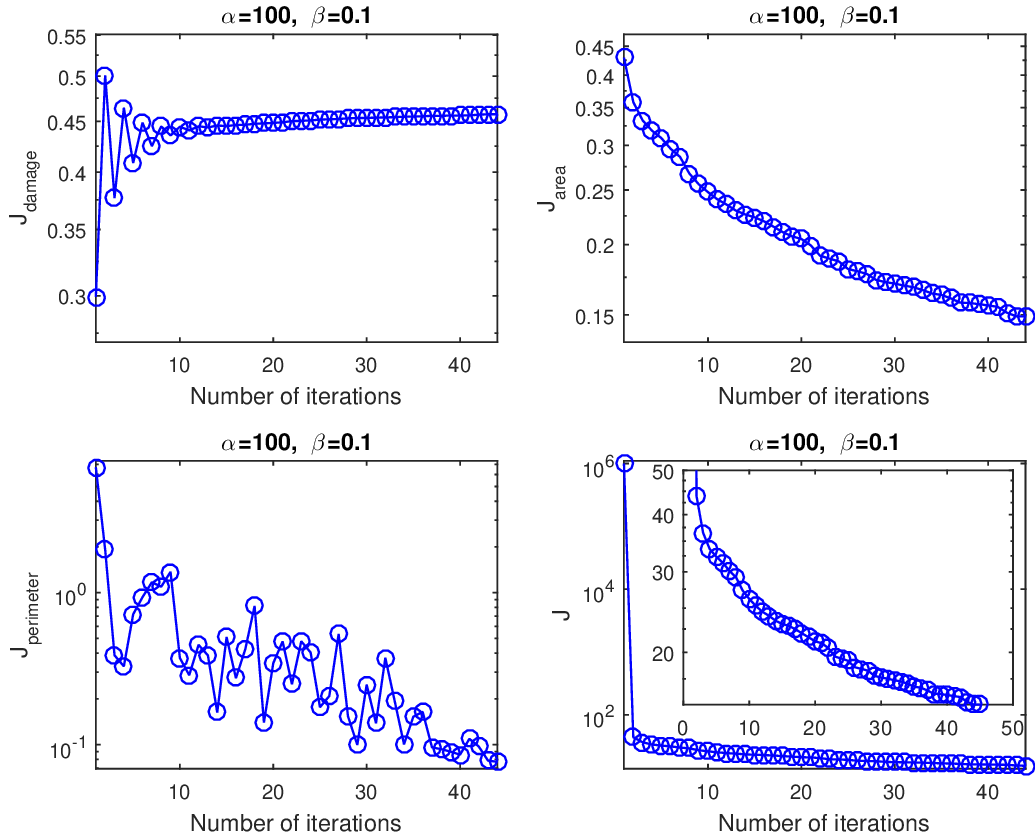}
\caption{Variation of the functionals $J_{damage}$, $J_{area}$, $J_{perimeter}$  and $J$ for $\alpha=100$ and $\beta=0.1$ (Experiment 1).}
\label{Fig_2}
\end{center}
\end{figure}

Figure~\ref{Fig_1} depicts the shape of the subdomain $\omega$ (plotted with light color) for $\alpha=100$ and $\beta=0.1$ at different iterations.
Variations of the functionals $J_{damage}$, $J_{area}$, $J_{perimeter}$ and $J$ for the weights $\alpha=100$ and $\beta=0.1$ are presented in Figure~\ref{Fig_2}. We note that $J_{damage}$ has an oscillatory behavior in the first 10 iterations of the optimization, and then tends to a stabilizing value. At the same time, the functional $J_{area}$ is continuously decreasing, whereas $J_{perimeter}$ although presents a general decreasing tendency has spurious fluctuations during the whole iterative process. These fluctuations being of small magnitude, they do not affect the decreasing evolution of the global functional $J$.

\begin{figure}[!tbhp]
\begin{center}
\includegraphics[scale=0.68]{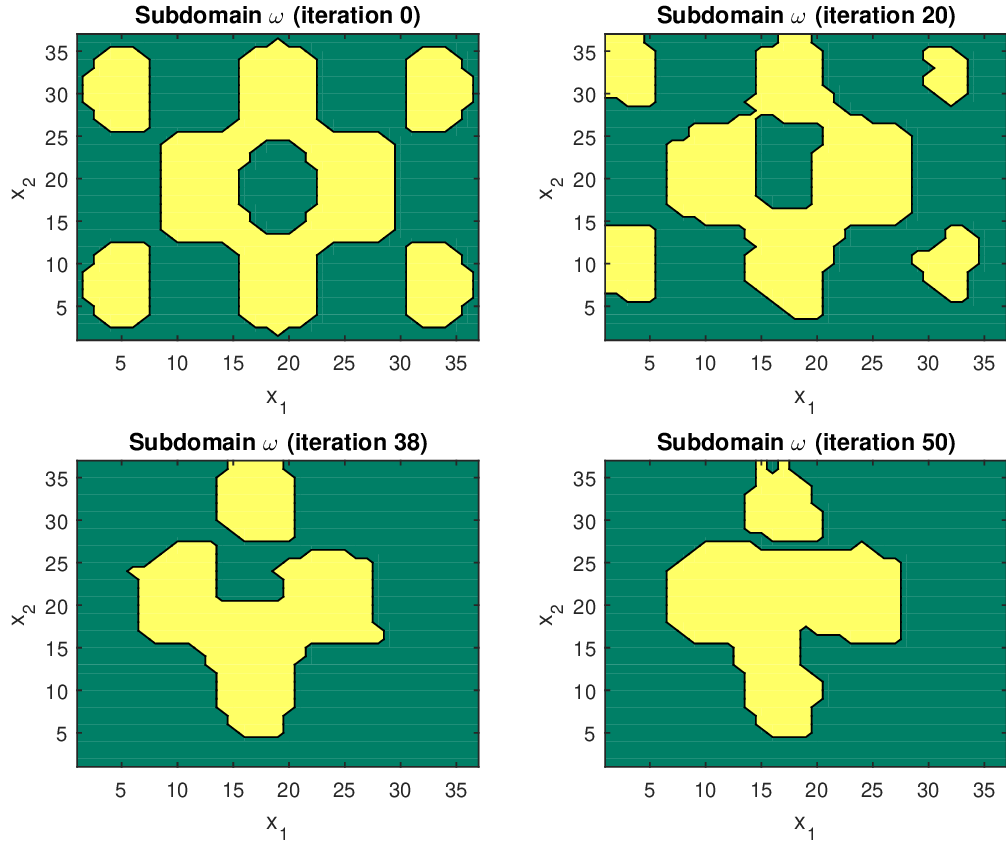}
\caption{Variation of the subdomain  $\omega$ (marked with light color) for $\alpha=0$ and $\beta=100$ at iterations $0, 20, 38$ and $50$   (Experiment 1).}
\label{Fig_3}

\bigskip

\bigskip

\includegraphics[scale=0.68]{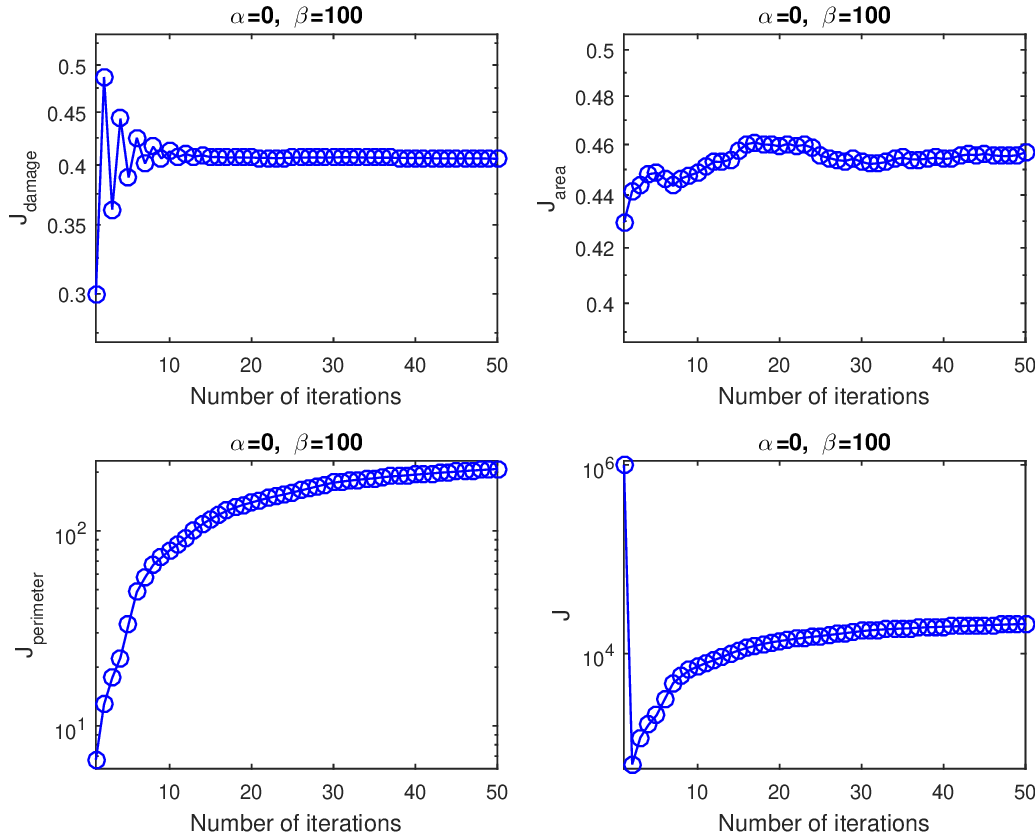}
\caption{Variation of the functionals $J_{damage}$, $J_{area}$, $J_{perimeter}$  and $J$ for $\alpha=0$ and $\beta=100$ (Experiment 1).}
\label{Fig_4}
\end{center}
\end{figure}

Analogous plots are obtained in Figures~\ref{Fig_3}--\ref{Fig_4} for $\alpha=0$, and $\beta=100$. In this case (with the new values for the weights, $\alpha$ and $\beta$) the algorithm is not convergent. The both functionals $J_{area}$ and $J_{perimeter}$ indicate an increasing tendency (even a strictly increasing for $J_{perimeter}$) that induce a similar bevaviour for $J$. All these functionals tend to stabilize to a certain value.

\bigskip

\textbf{Experiment $2$}.
In this experiment, we have used the same initial condition, $y_0(x_1,x_2)$, and the same initial shape function $\phi_0(x_1,x_2)$ as in Experiment 1, but we have chosen an asymmetric form for the function $\kappa(x,x') \equiv |\kappa_1(x_1,x_2) \kappa_2(x'_1,x'_2)|$ with
$$\kappa_1(x_1,x_2) = 500\sin(3\pi x_1)\cos(5\pi x_2)\exp(-(x_1 - x_2 - 0.2)^2 - 3(x_1 - x_2 - 0.8)^2)\,,\quad\mbox{and}$$
$$\kappa_2(x'_1,x'_2) = 500\sin(5\pi x'_1)\cos(3\pi x'_2)\exp(-5(x'_1 - 0.2)^2 - (x'_2 - 0.8)^2)\,.$$

\begin{figure}[!tbhp]
\begin{center}
\includegraphics[scale=0.68]{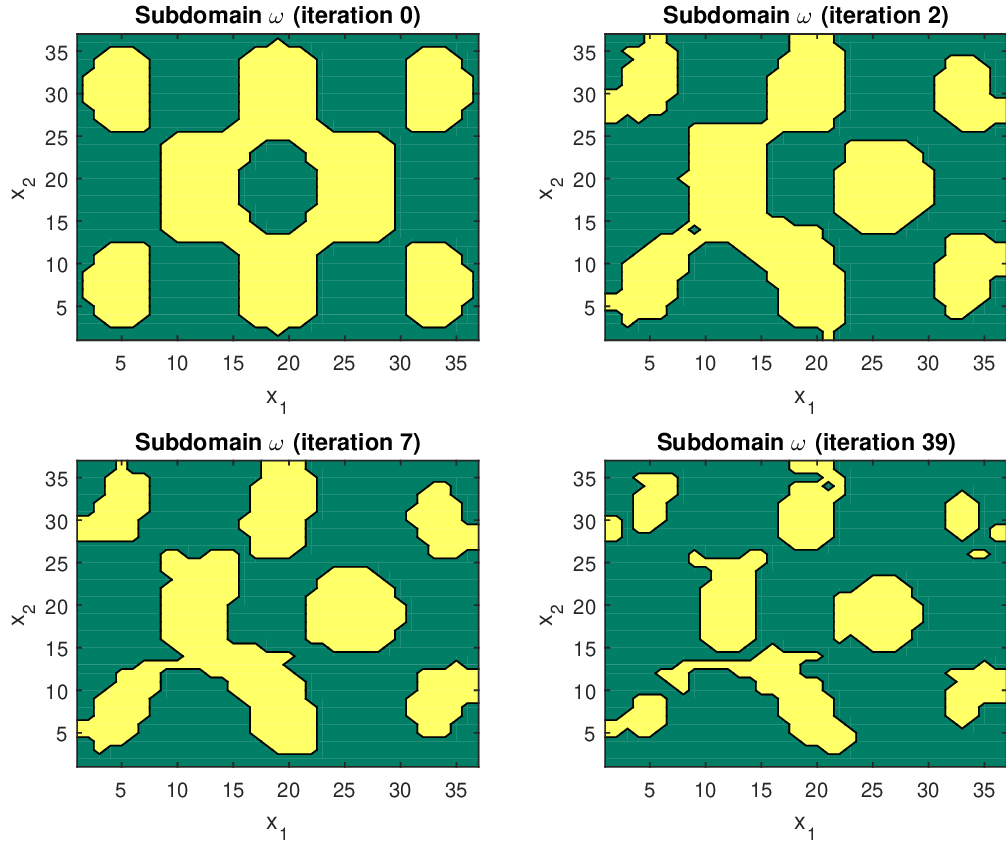}
\caption{Variation of the subdomain  $\omega$ (marked with light color) for $\alpha=100$ and $\beta=0.1$ at iterations $0, 2, 7$ and $39$  (Experiment 2).}
\label{Fig_5}

\bigskip

\bigskip

\includegraphics[scale=0.68]{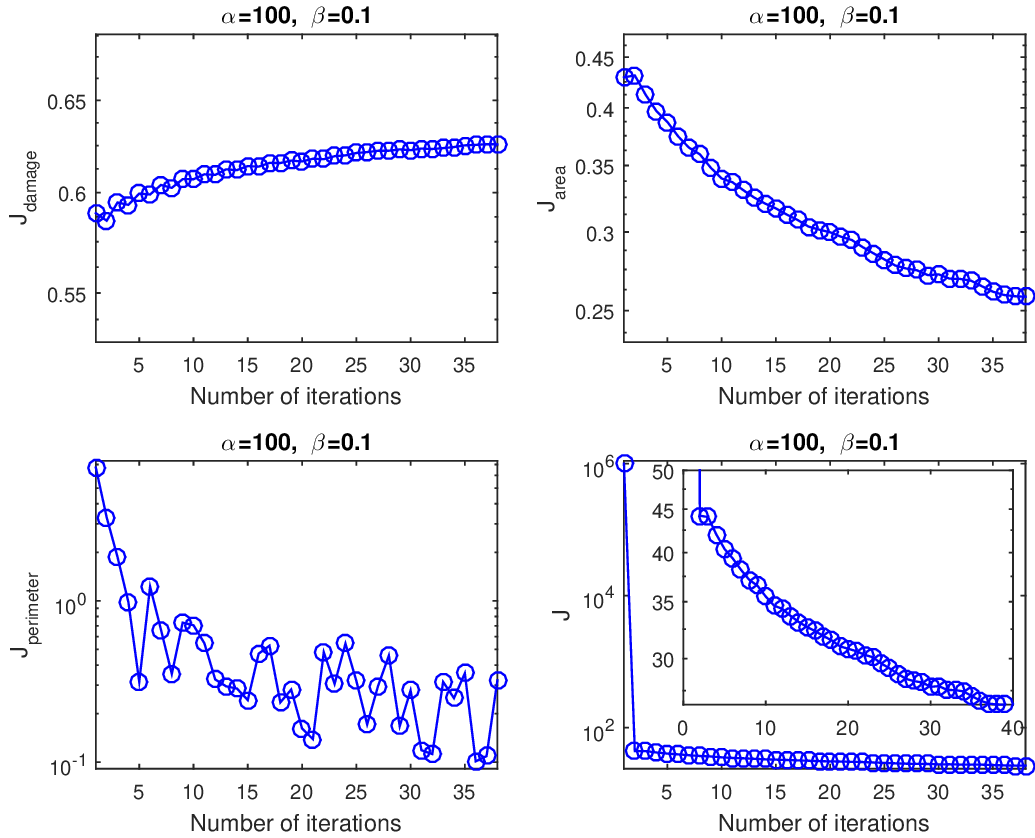}
\caption{Variation of the functionals $J_{damage}$, $J_{area}$, $J_{perimeter}$  and $J$ for $\alpha=100$ and $\beta=0.1$ (Experiment 2).}
\label{Fig_6}
\end{center}
\end{figure}

The shape of the subdomain $\omega$ (where the control acts, area marked with light color) for $\alpha=100$ and $\beta=0.1$ at different iterations is illustrated in Figure~\ref{Fig_5}.
The variation of the functionals $J_{damage}$, $J_{area}$, $J_{perimeter}$ and $J$ for $\alpha=100$ and $\beta=0.1$ is presented in Figure~\ref{Fig_6}. In this case the algorihm is convergent in 39 iterations. Zoomed areas in Figures~\ref{Fig_2} and \ref{Fig_6} (bottom--right) show more clear the decreasing of the global functional $J$.

\begin{figure}[!tbhp]
\begin{center}
\includegraphics[scale=0.90]{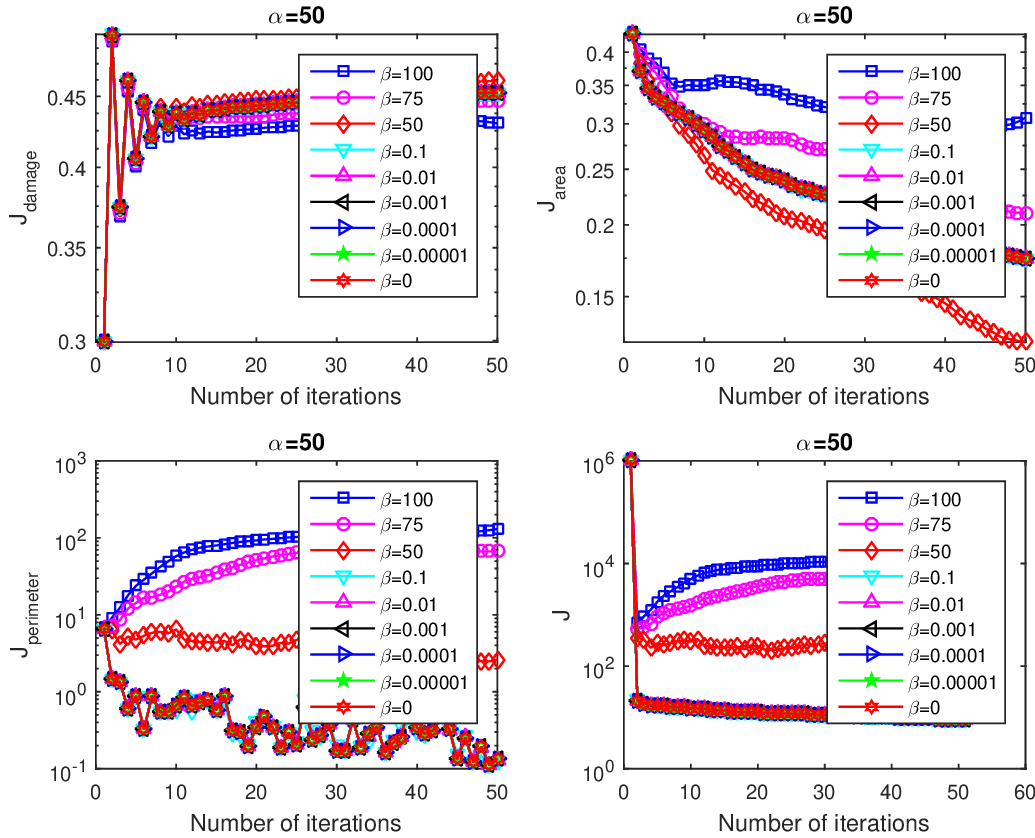}
\caption{Variation of the functionals $J_{damage}$, $J_{area}$, $J_{perimeter}$ and $J$ for $\alpha=50$ and several values of $\beta$ (Experiment 2).}
\label{Fig_7}
\end{center}
\end{figure}

\begin{figure}[!tbhp]
\begin{center}
\includegraphics[scale=0.90]{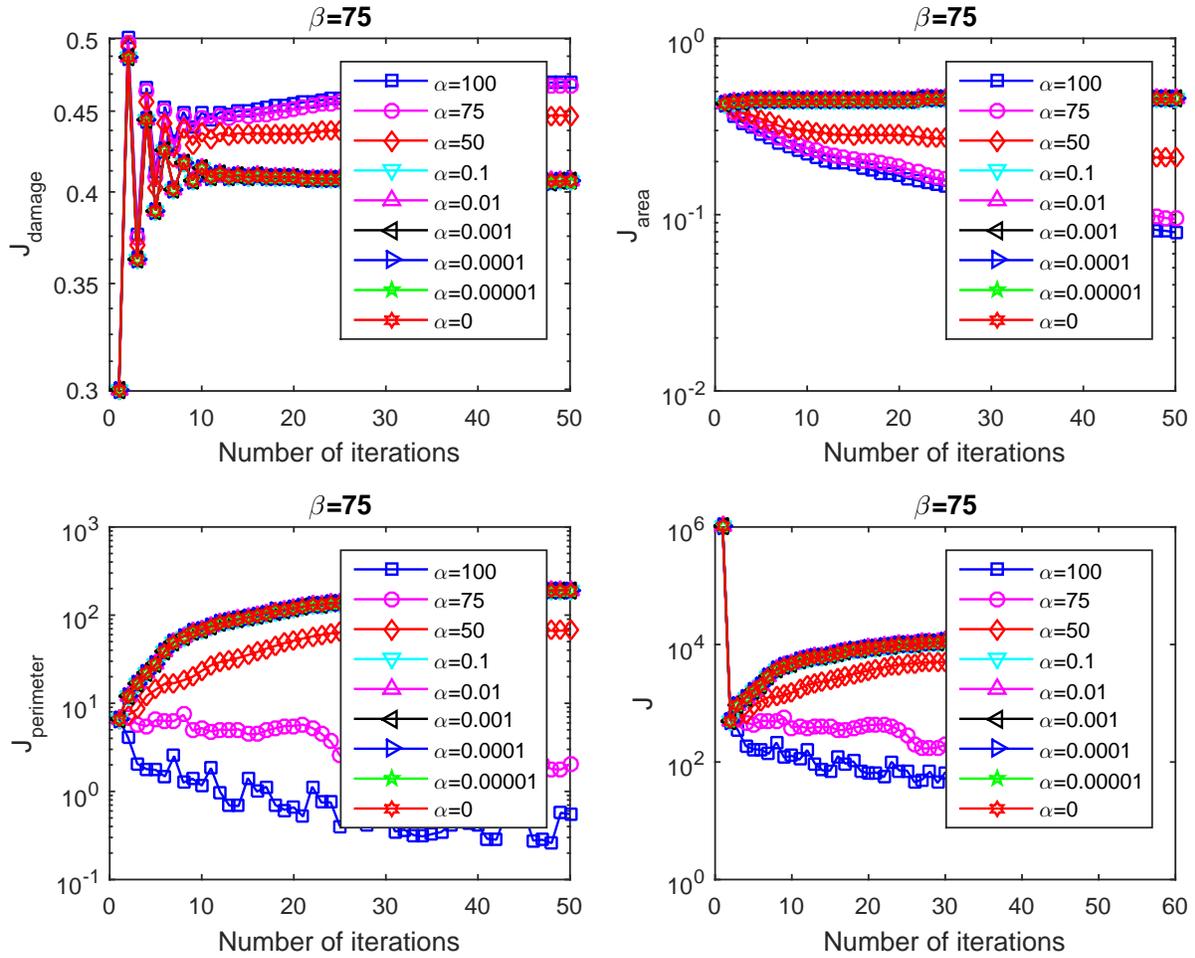}
\caption{Variation of the functionals $J_{damage}$, $J_{area}$, $J_{length}$ and $J$ for $\beta=75$ and several values of $\alpha$ (Experiment 2).}
\label{Fig_8}
\end{center}
\end{figure}

The last two figures (Figure~\ref{Fig_7} and Figure~\ref{Fig_8}) give an idea about the robustness of the minimization algorithm, when the weight parameters in the global functional $J$, $\alpha$ and $\beta$, vary in a certain range of values. Thus, Figure~\ref{Fig_7} shows the variations of the functionals $J_{damage}$, $J_{area}$, $J_{perimeter}$ and $J$ during iterative process, when  $\alpha=50$ and $\beta$ takes successively values in the set ${\cal{\cal W}}=\{0, 0.00001, 0.0001, 0.001$, $0.01$, $0.1, 50, 75, 100 \}$.
Analogously, Figure~\ref{Fig_8} presents the variations of the functionals $J_{damage}$, $J_{area}$, $J_{perimeter}$ and $J$ during iterative method, when  $\beta=75$ and $\alpha$ takes values in the set ${\cal{\cal W}}$.

\section{Final remarks}\label{AD_s5}

This work has proposed a novel regional control strategy of
minimizing the total cost of the damages produced by an alien
predator population and to reduce  this population. The dynamics
of the predators is described by a prey-predator system with local
or nonlocal reaction terms. A sufficient condition for the
zero-stabilizability (eradicability) of predators is given in
terms of the sign of the  principal eigenvalue of an appropriate
operator that is not self-adjoint, and a stabilizing feedback
control with a very simple structure is indicated. The
minimization related to such a feedback control is treated for a
closely related minimization problem viewed as a regional control
problem. The level set method has  been adopted  for handling the
Minkowski functionals of the relevant subregion.

An iterative algorithm to decrease the total cost was obtained and
numerical results showed the effectiveness of the theoretical
results. Several numerical simulations have been carried out for
the prey-predator system with nonlocal reaction terms. In this
case, the numerical response to predation reflects the
interactions among individuals in an actual habitat, when the
predators from position $x'$ that captured preys at position $x$,
will stay and produce offsprings at this new position (the
predators follow the prey). As a conclusion regarding the
numerical realization, the proposed  algorithm is strongly
affected by the values of the model parameters, and not in the
least, the results also depend on the resolution of the
discretisation (both space and time steps). Its convergence is
attained when an appropriate selection of the weight parameters,
$\alpha$ and $\beta,$ is done. These selections should maintain a
balance of the functionals $J_{damage}$, $J_{area}$ and
$J_{perimeter}$ with respect to their order of magnitude.

Finally, one may  remark that a spatially structured SIR problem
may be described by the same system, and the above mentioned
minimization problem may be viewed as the problem of minimizing
the effects of an epidemics by regional controls.


\section{Appendix}\label{AD_s6}
We establish here some auxiliary results. Consider the following eigenvalue problem
\begin{equation}\label{AD_app_eq_15}
\left\{
\begin{array}{ll}
-d\Delta\psi(x) + \eta(x)\psi(x) - K(x) (B\psi )(x) = \lambda\psi(x)\,, & x\in\Omega\,,\\
\partial_\nu\psi(x)=0\,, & x\in\partial\Omega\,.
\end{array}
\right.
\end{equation}
Here $d\in (0,+\infty)$, $\eta, K\in L^\infty(\Omega)$, $K(x)\ge 0$ a.e. in $\Omega$, and $B$ satisfies (A2).
\vspace{2mm}

\bf Lemma 1. \rm Problem (\ref{AD_app_eq_15}) has a simple eigenvalue $\lambda_1\in \mathbb{R}$, which corresponds to a positive eigenfunction. None of the other eigenvalues corresponds to a positive eigenfunction.
\vspace{2mm}

$\lambda_1$ is called the principal eigenvalue for (\ref{AD_app_eq_15}).
\vspace{2mm}

\bf Remark. \rm There exists $\psi_1$ an eigenfunction for (\ref{AD_app_eq_15}), corresponding to $\lambda_1$, such that
$$\psi_1(x)\ge\zeta > 0\,,\quad\forall x\in \Omega $$
(actually, $\psi_1\in C(\overline{\Omega}))$.
\vspace{2mm}

{\it Proof of Lemma.} Let us use the Krein-Rutman theorem. We may assume without loss of generality, that there exists
$\eta_0$ such that
$$\eta(x)\ge\eta_0>\|K\|_\infty \|B\|$$
a.e. $x\in \Omega $ (where $\| K\|_{\infty }=\| K\|_{L^{\infty }(\Omega )}, \ \|B\|=\|B\|_{L(L^2{\Omega}))}$). If this hypothesis is not satisfied, then we reduce our problem to this situation by translating $\lambda$.

Let $X=L^\infty(\Omega)$, $C=\{w\in L^\infty(\Omega)\,;\ w(x)\ge 0 \ \mbox{a.e.}\ x\in\Omega\}$, and
${\cal {\cal T}}:X\to X$ given by
$${\cal {\cal T}}f=\psi\,,$$
where $f\in L^\infty(\Omega)$ and $\psi$ is the unique solution to
\begin{equation}\label{AD_app_eq_16}
\left\{
\begin{array}{ll}
-d\Delta\psi(x) + \eta(x)\psi(x) - K(x)(B\psi )(x) = f(x)\,, & x\in\Omega\,,\\
\partial_\nu\psi(x) = 0\,, & x\in\partial\Omega\,.
\end{array}
\right.
\end{equation}
We have that $X$ is a Banach space and $C$ is a solid cone ($C$ is a closed convex cone with nonempty interior).

Let us prove that ${\cal {\cal T}}$ is a compact linear operator which is strictly positive (i.e., if $f\in C$ and
$f\ne 0_X$, then ${\cal {\cal T}}f\in\mbox{Int}\,C$). It is obvious that the hypotheses in Lax-Milgram lemma are satisfied if we view $f$ as an element of $L^2(\Omega)$. Hence, (\ref{AD_app_eq_16}) has a unique weak solution
$\psi\in W^{2,2}(\Omega)$.

Actually, by Theorem~IX.26 in \cite{brezis} we get that $\psi\in W^{2,2}(\Omega)$, and there exists a positive constant $\tilde{c}_2$ such that
$$\|\psi\|_{W^{2,2}(\Omega)} \le \tilde{c}_2\|f\|_{L^2(\Omega)}$$
($\tilde{c}_2$ is independent of $f$).

Finally, we get that there exists $\tilde{c}>0$ such that
$$\|{\cal {\cal T}}f\|_{W^{2,2}(\Omega)}\le \tilde{c}\|f\|_{L^\infty(\Omega)}\,.$$
Since $W^{2,2}(\Omega)\subset C(\overline{\Omega})$ continuously,
we obtain that ${\cal {\cal T}}f\in C(\overline{\Omega})\subset
L^\infty(\Omega)$, and ${\cal {\cal T}}$ is linear and bounded. On
the other hand, since the embedding $W^{2,2}(\Omega)\subset
C(\overline{\Omega })$ is compact (see \cite{brezis}), we may
infer that ${\cal {\cal T}}$ is a compact linear operator. Let us
prove that ${\cal {\cal T}}$ is strictly positive. Let $f\in
L^\infty(\Omega)$, $f(x)\ge 0$ a.e. $x\in\Omega$, and $f\ne 0_X$.
Let us prove that $\psi$, the weak solution to
(\ref{AD_app_eq_16}) satisfies $\psi(x)\ge 0$ a.e. $x\in \Omega$.
Indeed, if we consider $\psi^-\in W^{1,2}(\Omega)$ then
$$d\int_\Omega \nabla\psi(x)\cdot\nabla\psi^-(x)\,dx + \int_\Omega \eta(x)\psi(x)\psi^-(x)\,dx
-\int_\Omega K(x)(B\psi )(x)\psi^-(x)\,dx
=\int_\Omega f(x)\psi^-(x)\,dx \ge 0\,.$$
This implies that
$$-d\int_\Omega|\nabla\psi^-(x)|^2\,dx - \int_\Omega \eta(x)|\psi^-(x)|^2\,dx
+ \int_\Omega K(x)(B\psi^- )(x)\psi^-(x)\,dx\ge 0\,,$$
and consequently,
$$\|K\|_{\infty }\cdot \|B\|\cdot \|\psi^-\|^2_{L^2(\Omega )}\geq
\int_\Omega \eta(x)|\psi^-(x)|^2\,dx\ge \eta_0 \|\psi^-\|^2_{L^2(\Omega)}\,.$$
Since $\eta_0>\|K\|_\infty\|B\|$, it follows that
$\psi^-=0_{L^2(\Omega)}$, and so $\psi(x)\ge 0$ a.e. $x\in\Omega$.

Let us prove that $\psi(x)>0$,\    for any $
x\in\overline{\Omega}$ (recall that $\psi\in
C(\overline{\Omega}))$. Indeed, since $K(x)(B\psi(\cdot))(x)\geq
0$ a.e. $x\in\Omega$, we conclude that
\begin{equation*}
\left\{
\begin{array}{ll}
-d\Delta \psi(x) + \|\eta\|_\infty\psi(x)\ge 0\,, & x\in\Omega\,,\\
\partial_\nu\psi(x) = 0\,, & x\in\partial\Omega\,.
\end{array}
\right.
\end{equation*}
It follows that for any $q\in [1,+\infty)$, there exists $c_q>0$ such that
$$\|\psi\|_{L^q(\Omega )} \leq c_q\cdot \mbox{inf}_{x\in\Omega}\, \psi (x)$$
(see Lemma 2.1 in \cite{cheny}).

If $\mbox{inf}_{x\in\Omega}\,\psi (x)=\mbox{min}_{x\in\overline{\Omega}}\,\psi (x)>0$, then we get that
${\cal {\cal T}}f > 0$.

Indeed, if we assume, by contradiction, that
$\mbox{inf}_{x\in\Omega}\,\psi (x)= 0$, then it follows that
$\|\psi\|_{L^q(\Omega )} = 0$, for any $q\in [1,+\infty)$. We get
that $\psi(x)=0$ a.e. $x\in\Omega$, and so $\psi(x)=0$, $\forall\,
x\in\Omega$. This implies that $f=0_X$  which is a  contradiction.
Hence, ${\cal {\cal T}}$ is strictly positive.

\medskip

By Theorem 1.2 in \cite{du} we get that the spectral radius of ${\cal {\cal T}}$ satisfies $r({\cal {\cal T}}) > 0$, and
$r({\cal {\cal T}})$ is a simple eigenvalue with an eigenvector $\psi\in\mbox{Int}\, C$; there is no other eigenvalue with positive eigenvector.
The conclusion of Lemma 1 is now obvious.

\medskip

The second result concerns the principal eigenvalues for
\begin{equation}\label{AD_app_eq_17}
\left\{
\begin{array}{ll}
-d\Delta \psi(x) + \eta(x)\psi(x) - (K(x)+\varepsilon)(B\psi )(x) = \lambda\psi(x)\,, & x\in\Omega\,,\\
\partial_\nu\psi(x) = 0\,, & x\in\partial\Omega\,,
\end{array}
\right.
\end{equation}
where $\varepsilon > 0$. We denote by $\lambda_1(\varepsilon)$ the principal eigenvalue to (\ref{AD_app_eq_17}).
\vspace{2mm}

\bf Lemma 2. \rm $$\lim_{\varepsilon\to 0} \lambda_1(\varepsilon) = \lambda_1\,.$$
\vspace{2mm}

\it Proof. \rm Let $\psi_\varepsilon$ be the positive eigenfunction corresponding to $\lambda_1(\varepsilon)$, and
satisfying
$$\|\psi_\varepsilon\|_{L^2(\Omega)} = 1\,.$$
Let us prove that if $0\le \varepsilon_1 < \varepsilon_2$, then
$$\lambda_1(\varepsilon_1)\ge \lambda_2(\varepsilon_2)\,.$$

Consider $\tilde{\lambda}_1(\varepsilon_1)$ the principal eigenvalue for
\begin{equation}\label{AD_app_eq_18}
\left\{
\begin{array}{ll}
-d\Delta \psi(x) + \eta(x)\psi(x)  + (B^*((K(\cdot )+\varepsilon_1)\psi (\cdot )))(x) = \lambda\psi (x)\,, & x\in \Omega\,,\\
\partial_\nu\psi(x) = 0\,, & x\in\partial\Omega\,.
\end{array}
\right.
\end{equation}
The existence and basic properties related to it follow as for
(\ref{AD_app_eq_15}). In fact
$\lambda_1(\varepsilon)\in\mathbb{R}$, and let
$\tilde{\psi}_{\varepsilon_1}$ be  a corresponding and positive
eigenfunction.

Using (\ref{AD_app_eq_17}) and (\ref{AD_app_eq_18}) we get that
\begin{equation}\label{AD_app_eq_19}
-(\varepsilon_2 - \varepsilon_1)\int_\Omega (B\psi_{\varepsilon_2})(x)\cdot
\tilde{\psi}_{\varepsilon_1}(x)\,dx
= (\lambda_1(\varepsilon_2) - \tilde{\lambda}_1(\varepsilon_1))
\int_\Omega\psi_{\varepsilon_2}(x) \tilde{\psi}_{\varepsilon_1}(x)\,dx\,,
\end{equation}
and
$$0= (\lambda_1(\varepsilon_1) - \tilde{\lambda}_1(\varepsilon_1))\int_\Omega \psi_{\varepsilon_1}(x)\tilde{\psi}_{\varepsilon_1}(x)\,dx\,.$$
Since $\int_\Omega \psi_{\varepsilon_1}(x)\tilde{\psi}_{\varepsilon_1}(x)\,dx > 0$, we may conclude that
$\tilde{\lambda}_1(\varepsilon_1) = \lambda_1(\varepsilon_1)$.

On the other hand, since $\int_\Omega(B\psi_{\varepsilon_2})(x)\tilde{\psi}_{\varepsilon_1}(x)\,dx \geq 0$,
and $\int_\Omega \psi_{\varepsilon_2}(x)\tilde{\psi}_{\varepsilon_1}(x)\,dx >0$, we get by (\ref{AD_app_eq_19}) that
$\lambda_1(\varepsilon_2)\le \lambda_1(\varepsilon_2)$.

We also conclude that $\lambda_1(\varepsilon)\le \lambda_1$\,, for any $\varepsilon > 0$. We may infer that there exists
$\lim_{\varepsilon\to 0} \lambda_1(\varepsilon) = \tilde{\lambda}_1 \le \lambda_1$.

Let us prove that actually we have equality. If $0<\varepsilon<1$, then by (\ref{AD_app_eq_17}) we get that
$$d\int_\Omega|\nabla\psi_\varepsilon|^2\,dx + \int_\Omega \eta(x)|\psi_\varepsilon|^2\,dx\le \lambda_1(\varepsilon)
+(\|K\|_{\infty }+\varepsilon)\cdot\|B\|\,,$$
and consequently, $\psi_\varepsilon$ is bounded in $W^{1,2}(\Omega)$. Therefore, there exists a sequence
$(\psi_{\varepsilon_n})$\ $(\varepsilon_n\to 0)$, such that
$$\psi_{\varepsilon_n}\rightharpoonup \psi_0\quad\mbox{in}\ W^{1,2}(\Omega)\,,\quad\mbox{and}\quad
\psi_{\varepsilon_n}\to \psi_0\quad\mbox{in}\ L^2(\Omega)\,,$$
which implies that $\|\psi_0\|_{L^2(\Omega)}=1$, and $\psi_0(x)\ge 0$ a.e. $x\in\Omega$. Since
$\psi_{\varepsilon_n}$ satisfies
\begin{eqnarray*}
& & d\int_\Omega \nabla\psi_{\varepsilon_n}(x)\cdot \nabla\psi(x)\,dx + \int_\Omega \eta(x)\psi_{\varepsilon_n}(x)\psi(x)\,dx\\
&-& \int_\Omega(K(x) + \varepsilon_n)(B\psi_{\varepsilon_n})(x)\psi(x)\,dx
=\lambda_1(\varepsilon)\int_\Omega\psi_{\varepsilon_n}(x)\psi(x)\,dx\,,
\end{eqnarray*}
for any $\psi \in W^{1,2}(\Omega )$, we may pass to the limit, and obtain that $\psi_0$ is a weak solution to (\ref{AD_app_eq_15}) corresponding to
$\tilde{\lambda}_1$, i.e. $\psi_0$ is a nonnegative eigenfunction for (\ref{AD_app_eq_15}), corresponding to
$\lambda : = \tilde{\lambda}_1$, and so $\tilde{\lambda}_1$ is an eigenvalue. By Lemma 1 we conclude that $\tilde{\lambda}_1 = \lambda_1$.

\medskip

Using the same method we have used for Lemma 2  may prove that
\vspace{2mm}

\bf Lemma 3. \rm The mapping $\gamma\mapsto \lambda^\omega_{1\gamma}$ is strictly increasing.


\bigskip

\section*{Acknowledgements}
This work has been carried out in the framework of the   European  COST
project CA16227: ``Investigation and Mathematical Analysis of Avant-garde Disease Control via Mosquito Nano-Tech-Repellents''.

\pagebreak


\end{document}